\documentclass{article}
\usepackage[left=2cm,top=2cm,right=3cm,bottom=1cm, footskip=.5cm]{geometry}
\usepackage{amsmath}
\usepackage{graphicx}
\usepackage{tikz}
\usepackage{fancyvrb}
\usepackage{pdfpages}
\usepackage[makeroom]{cancel}
\usepackage{amssymb} 
\usepackage[parfill]{parskip}
\usepackage{mathtools}
\usepackage{hyperref}
\usepackage{tikz-cd}
\usepackage{setspace}
\usepackage{amsthm}
\usepackage{color}
\usetikzlibrary{%
	matrix,%
	calc,%
	arrows%
}

\usepackage[style = numeric, backend = bibtex]{biblatex}
\newtheoremstyle{colon}%
{}
{}
{\itshape}%bodyfont
{}%indent
{\bfseries}%headfont
{:}%head punctuation
{ }%space after head
{}

\theoremstyle{colon}

\newtheorem{Lemma}{Lemma}[subsection]
\newtheorem{Theorem}[Lemma]{Theorem}
\newtheorem{Proposition}[Lemma]{Proposition}

\title{Completion preserves homotopy fibre squares of connected nilpotent spaces}

\begin{document}
	
	\date{}
	\author{A. Ronan}
	\maketitle
	
\begin{abstract}
We prove that completion preserves homotopy fibre squares of connected nilpotent spaces. As a consequence, we deduce the Hasse fracture square associated to a connected nilpotent space. Along the way, we give a quick proof of the well-known result that if the base of a fibre sequence has the homotopy type of a CW complex, then the total space has the homotopy type of a CW complex iff each fibre does.
\end{abstract}

\section{Introduction}

Let T be a non-empty set of primes. By a functorial completion, we mean a pair $(F,\alpha)$ where F is a functor that takes nilpotent spaces to nilpotent spaces, and $\alpha: 1 \to F$ is a natural transformation such that, for all $X$, $\alpha_X: X \to F(X)$ is a T-completion of $X$. For example, such a completion can be obtained via Bousfield localisation at the homology theory $H_{*}(-;\bigoplus_{p \in T} \mathbb{F}_p)$, or by using the Bousfield-Kan completion functor. By a homotopy fibre square we mean a strictly commutative square such that the canonical map to the double mapping path space, which we denote by $N(f,g)$, is a weak equivalence.  Our main result now states:

\begin{Theorem} Let $X \to \hat{X}_T$ denote a functorial completion. If a commutative square of connected nilpotent spaces:
	
	\[ \begin{tikzcd}
	P \arrow[swap]{d}{k} \arrow{r}{h} &  Y \arrow{d}{g} \\
	X \arrow[swap]{r}{f} & A\\
	\end{tikzcd} \]
	
	is a homotopy fibre square, then so is:
	
	\[ \begin{tikzcd}
	\hat{P}_T \arrow[swap]{d}{\hat{k}_T} \arrow{r}{\hat{h}_T} &  \hat{Y}_T \arrow{d}{\hat{g}_T} \\
	\hat{X}_T \arrow[swap]{r}{\hat{f}_T} & \hat{A}_T \\
	\end{tikzcd} \]	
\end{Theorem}

Similar results to Theorem 1.0.1 can be found in the literature. In particular, in \cite[Theorem 1.1]{ref1}, Farjoun proves an analogue of Theorem 1.0.1 in the case of disconnected spaces. However, one can only conclude directly from Theorem 1.1 of \cite{ref1} that the comparison map $\hat{P}_T \to N(\hat{f}_T, \hat{g}_T)$ has homotopically discrete fibre. Nevertheless, it seems likely that only a small modification of the proof given in \cite{ref1} would be enough to derive our Theorem 1.0.1 in the case where the functorial completion in question is Bousfield-Kan's $\mathcal{R}_{\infty}$ functor. On the other hand, one bonus of our proof is that it holds for any functorial completion, such as the one obtained via Bousfield localisation.

Our main reason for being interested in Theorem 1.0.1 is that it provides a natural context in which to prove the following well-known fracture theorem, sometimes known as the Hasse square:

\begin{Theorem}
	Let $X$ be a $T$-local connected nilpotent space. Then any commutative square:
	
	\[ \begin{tikzcd}
	X \arrow{r}{\hat{\phi}} \arrow[swap]{d}{\psi} & \hat{X}_T \arrow{d}{\phi} \\
	X_0 \arrow[swap]{r}{\hat{\phi}_0} & (\hat{X}_T)_0 \\
	\end{tikzcd}
	\]
	
	with $\hat{\phi}$ a T-completion and $\psi,\phi$ rationalisations, is a homotopy fibre square.
	
\end{Theorem}

For a standard proof, see \cite[Theorem 4.4]{ref7}. We explain how the Hasse square can be deduced from Theorem 1.0.1 at the end of this paper.

In general, we like to think of the T-completion of a nilpotent space $X$ as being constructed by first using chosen central series for the homotopy groups of $X$ to replace $X$ by a Postnikov tower. Then, a T-completion of $X$ can be defined by induction up the tower, using T-completions of Eilenberg-MacLane spaces. This is the approach taken in May and Ponto's book \cite{ref2}. This construction is functorial only up to homotopy, and this introduces a number of difficulties in trying to prove Theorem 1.0.1. Our solution to these problems is to take the coattaching maps of Postnikov towers to be cofibrations, which allows us to make more squares commute when defining a T-completion of $X$. It turns out that this introduces sufficient functoriality to be able to prove Theorem 1.0.1.

Along the way, we will review some results about well-pointed spaces and spaces with the (unbased) homotopy type of a CW complex. In particular, we give the most direct proof we know of the following well-known, see \cite{ref8}, result:

\begin{Theorem}
	Let $B$ be a connected space with the homotopy type of a CW complex and let $p: E \to B$ be an h-fibration with fibre $F$. Then $E$ has the homotopy type of a CW complex iff $F$ has the homotopy type of a CW complex.
\end{Theorem}

One direction of our proof can be viewed as an application of the mixed model structure of M.Cole, \cite{ref5}. We will also discuss the actions of fundamental groups on a composite of fibrations and some basic results about the category $\mathcal{B}_T$ of T-complete abelian groups, as well as their nilpotent analogues.

\textbf{Counterexamples}

We now say a few words about the hypotheses of the theorem. Consideration of the path-space fibration associated to $K(G,1)$ for any $G$ with $\pi_2(\widehat{K(G,1)}_T) \neq 0$, shows that some connectivity assumptions on $P$ are necessary. Consideration of the two commutative squares:

\[ \begin{tikzcd}
S^1 \arrow{r}{i} \arrow[swap]{d}{i} & D^2 \arrow{d}{i} & & & S^1 \arrow{r}{i} \arrow[swap]{d}{i} & D^2 \arrow{d}{*} \\
D^2 \arrow[swap]{r}{i} & S^2 & & & D^2 \arrow[swap]{r}{*} & S^2 \\
\end{tikzcd} \]

shows that only one of them is a homotopy cofibre square, despite all component maps being homotopic. This explains why some degree of functoriality is also required.

Finally, we state a counterexample, due to Sullivan, when the spaces involved are not nilpotent. We assume that T is a single prime $p$, and that the functorial completion in question is obtained via Bousfield localisation of $\mathbb{F}_p$-homology. In fact, the spaces we consider are all '$\mathbb{F}_p$-good', so this is also a counterexample for Bousfield-Kan's $\mathcal{R}_{\infty}$ functor \cite[Ch. VII, Prop. 2.3, 3.6]{ref9}. The counterexample is based on a non-nilpotent space $Z$ satisfying the following properties, see \cite[pg. 104]{ref10} and \cite[Ch. VII. 3.6]{ref9} for details of the construction:

i) $\pi_1(Z) = \frac{\mathbb{Z}}{n \mathbb{Z}}$, where $n$ can be any integer dividing $p-1$, 

ii) $\pi_2(Z) = \hat{\mathbb{Z}}_p$, 

iii) $\pi_i(Z) = 0$ for $i \geq 3$, 

iv) $\Omega \hat{Z}_p \simeq \hat{S}^{2n-1}_p$. 

For large values of $n$ and $p$, it is clear from these properties that $p$-completion cannot preserve the fibre sequence:

\[ K(\hat{\mathbb{Z}}_p,2) \to Z \to K(\frac{\mathbb{Z}}{n \mathbb{Z}},1)
\]

\textbf{Notation}

We use throughout the notations and conventions of \cite{ref2}. Indeed, anybody who has read \cite{ref2} has more than enough background to understand this paper. For example, we have the definitions $\mathbb{E}_T G := \pi_1(\widehat{K(G,1)}_T), \mathbb{H}_T G := \pi_2(\widehat{K(G,1)}_T)$. If $G$ is abelian, then $\mathbb{H}_T G$ and $\mathbb{E}_T G$ are the first and zeroth derived functors of T-adic completion, respectively. In the abelian case, $\mathbb{H}_T G = Hom(\mathbb{Z}[T^{-1}] / \mathbb{Z}, G)$ and $\mathbb{E}_T G = Ext(\mathbb{Z}[T^{-1}] / \mathbb{Z}, G) $, which justifies the notation.

\newpage

If $\mathcal{C}$ is a class of abelian groups, then we call $G$ $\mathcal{C}$-nilpotent if $G$ has a central series such that the abelian subquotients are in $\mathcal{C}$. If $R$ is a ring, then we say that $G$ is $R$-nilpotent if it has a central series such that the abelian subquotients can be given the structure of an $R$-module. We say that $G$ is $fR$-nilpotent if we can give the abelian subquotients the structure of a finitely generated $R$-module. We sometimes use the term $fR$-module as shorthand for a finitely generated $R$-module.

We are interested in three model structures on the category of CGWH spaces, namely the Quillen, Hurewicz and mixed model structures \cite{ref5}. If we wish to be precise about which model structure a fibration belongs to, then we use the terms $q$-fibration, $h$-fibration and $m$-fibration, respectively. 

\newpage

\section{Preliminaries}

\subsection{Well-pointed spaces with the homotopy type of a CW complex}

We begin by reviewing the results which guarantee that all spaces we consider are well-pointed and have the homotopy type of a CW complex. We work throughout in the category of CGWH spaces and the next few results follow from the characterisation of an $h$-cofibration as an NDR-pair. 

\begin{Lemma}
	The pullback of an h-cofibration along an h-fibration is an h-cofibration.

\end{Lemma}

\textbf{Proof:} See \cite[Lemma 1.3.1]{ref2}. $\hfill \square$ 

The next two results are \cite[Lemmas 5 and 6]{ref3}.

\begin{Lemma}
	If $i:A \to B$ and $j: B \to C$ are maps and both $j$ and $k = ji$ are $h$-cofibrations, then $i$ is an $h$-cofibration.
\end{Lemma}

\textbf{Proof:} Firstly, $i$ is an inclusion by verification of the universal property. Let $(H,\lambda)$ and $(K, \mu)$ represent $(C,B)$ and $(C,A)$ as NDR-pairs. Define $\mu^{'}(c) =  \mu(c) + \sup_{t \in I} \lambda(K(c,t))$ . Define:

\[ L(b,t) = \begin{cases}
H(K(b,t),1) & \mu^{'}(b) \leq \frac{1}{2} \\
H(K(b, 2(1-\mu^{'}(b))t),1) & \mu^{'}(b) \geq \frac{1}{2} \\
\end{cases}
\]

Then $(L, 2\mu^{'})$ represents $(B,A)$ as an NDR-pair. $\hfill \square$

\begin{Lemma}
	If $E,B,X$ are well-pointed and $p$ is an $h$-fibration in the following pullback square:
	
	\[
	\begin{tikzcd}
	P \arrow{r} \arrow{d} & E \arrow[two heads]{d}{p} \\
	X \arrow{r} & B \\
	\end{tikzcd}
	\]
	
Then $P$ is well-pointed.
\end{Lemma}

\textbf{Proof:} This is a corollary of Lemma 2.1.1 and Lemma 2.1.2. $\hfill \square$

Next we turn our attention to spaces with the unbased homotopy type of a CW complex. Recall that there is an m-model structure on spaces where the weak equivalences are the q-weak equivalences and the fibrations are the $h$-fibrations,  see \cite{ref5} and \cite[Section 17.3]{ref2}. The $m$-cofibrant objects are precisely the spaces with the unbased homotopy type of a CW complex. We will give the shortest proof we know of the following well-known, \cite{ref8}, theorem:

\begin{Theorem} Let $p: E \to B$ be an $h$-fibration with fibre $F$ and connected base $B$. If $B$ is $m$-cofibrant, then $E$ is $m$-cofibrant iff $F$ is $m$-cofibrant.
\end{Theorem}

Since the h-model structure is proper, we can and will assume that $B$ is already a CW-complex. The m-model structure will play a simplifying role in the proof of Lemma 2.1.6 below. However, first we prove:

\begin{Lemma} If $E$ is $m$-cofibrant, then so is $F$.
\end{Lemma}

\textbf{Proof:} Consider the square:

\[ \begin{tikzcd}
\vert Sing(E) \vert \arrow[swap, two heads]{d}{|Sing(p)|} \arrow{r}{\simeq} & E \arrow[two heads]{d}{p} \\
\vert Sing(B) \vert \arrow[swap]{r}{\simeq} & B \\
\end{tikzcd}
\]

Since $p$ is an $h$-fibration it is a $q$-fibration, so $Sing(p)$ is a Kan fibration and $|Sing(p)|$ is an $h$-fibration, by \cite[Theorem 17.5.7]{ref2} or \cite[Theorem 4.5.25]{ref4}. Since the $h$-model structure is proper, there is an induced homotopy equivalence between $F$ and $|Sing(F)|$, which is a $CW$ complex. $\hfill \square$

\begin{Lemma} If $F$ is $m$-cofibrant, then so is $E$.
\end{Lemma}

\textbf{Proof:} Suppose that the $n$-skeleton of $B$ is defined by attaching maps $\alpha_i: \partial \Delta^n \to B_{(n-1)}$. Let $e_i$ denote the inclusion of the cell into $B_{(n)}$, $e_i: \Delta^n \to B_{(n)}$. Define $P_i^{\partial \Delta^n}$ to be the pullback of $\alpha_i$ along $p$, and $P_i^{\Delta^n}$ to be the pullback of $e_i$ along $p$. Then $P_i^{\partial \Delta^n} \simeq \partial \Delta^n \times F$ and $P_i^{\Delta^n} \simeq \Delta^n \times F$ and so both are $m$-cofibrant. We have the pushout:

\[ \begin{tikzcd}
\sqcup_i P_i^{\partial \Delta^n} \arrow{d} \arrow{r} & p^{-1}(B_{(n-1)}) \arrow{d} \\ \sqcup_i P_i^{\Delta^n} \arrow{r} & p^{-1}(B_{(n)}) \\
\end{tikzcd} \]

Since the pullback of an $h$-cofibration along an $h$-fibration is an $h$-cofibration, $P_i^{\partial \Delta^n} \to P_i^{\Delta^n}$ is an $h$-cofibration. Therefore, the LHS map is an $h$-cofibration between $m$-cofibrant objects and, hence, an $m$-cofibration by \cite[Prop. 17.3.4]{ref2}. Therefore, the RHS is also an $m$-cofibration, and so, inductively, $p^{-1}(B) = E$ is $m$-cofibrant.  $\hfill \square$

\subsection{Preliminary lemmas}

By the results of the previous subsection, if $f: X \to Y$ is a map between well-pointed spaces, then $Ff$ is well-pointed and we can define an action of $\pi_1(X)$ on $\pi_n(Ff)$ in the usual way. For more general maps $f: X \to Y$ we can define an action of $\pi_1(X)$ on $\pi_n(Ff)$ by using a Reedy cofibrant approximation to $f$. 

\newpage

We will need the following version of the relative Hurewicz theorem, which can be proved in the same way as other versions:

\begin{Lemma}
	Let $n \geq 1$. If $f:X \to Y$ is a map such that $Ff$ is $(n-1)$-connected and $\pi_1(X)$ acts trivially on $\pi_n(Ff)$, then $\eta: Ff \to \Omega Cf$ is an $n$-equivalence which induces an isomorphism on $\pi_n$.
\end{Lemma}

In the next section, we require that the coattaching maps of Postnikov towers are cofibrations, and to ensure this we use the following lemma:

\begin{Lemma}
	Let $f:X \to Y$ be a map between well-pointed spaces of the homotopy type of a CW complex such that $Ff \simeq K(A,n)$ and  $\pi_1(X)$ acts trivially on $A$. Then there exists an equivalence $X \to Fk$ over $Y$, for some cofibration $k: Y \to K(A,n+1)$. 
\end{Lemma}	
	
\textbf{Proof:}	We know that $\eta: Ff \to \Omega Cf$ induces an isomorphism on $\pi_i$ for $i \leq n$, by Lemma 2.2.1. Therefore, there is a cofibration $j: Cf \to K(A,n+1)$ which induces an isomorphism on $\pi_{n+1}$. We now consider the diagram:

\[\begin{tikzcd}
X \arrow{r}{\nu} \arrow[swap]{dr}{f} & Fi \arrow{r}{w} \arrow{d} & Fk \arrow{dl} \\
& Y \arrow[tail, swap]{d}{i} \arrow[tail]{dr}{k} & \\
& Cf \arrow[tail, swap]{r}{j} & K(A,n+1) \\ 
\end{tikzcd} \]
	
where $w$ is induced by $j$ and $\nu$ is the canonical map which induces $\eta$ on fibres. By construction, $w \nu$ is a weak equivalence, as desired. $\hfill \square$

We will also need the following result which concerns the action of fundamental groups on the fibres of a composite of fibrations:

\begin{Lemma}
	Consider a triangle of fibrations with well-pointed fibres:
	
\[	\begin{tikzcd}
	X \arrow[swap, two heads]{d}{f} \arrow[two heads]{dr}{h} \\
	Y \arrow[swap, two heads]{r}{g} & Z \\
	\end{tikzcd} \]

Let $F_1,F_2,F_3$ denote the fibres of $f,h$ and $g$ respectively and note that $F_1$ is the fibre of $F_2 \twoheadrightarrow F_3$. Let $\gamma \in \pi_1(X)$. Then there is a commutative square:

\[ \begin{tikzcd}
F_2 \arrow{r} \arrow[swap]{d}{\gamma_2} & F_3 \arrow{d}{\gamma_3} \\
F_2 \arrow{r} & F_3 \\
\end{tikzcd} \]

such that $\gamma_2, \gamma_3$ represent $\gamma$ and $f_{*}(\gamma)$ respectively, and the induced map, $\gamma_1$, between fibres represents $\gamma$.

\end{Lemma}

\textbf{Proof:} Let $\gamma$ also denote a loop in $X$ representing $\gamma$. We first define $\gamma_3$ in the usual way by constructing a homotopy $H: F_3 \times I \to Y$. We then define $\gamma_2$ by constructing a homotopy $G$ as in the following lifting problem:

\[ \begin{tikzcd}
F_2 \times \{0\} \cup * \times I  \arrow{r}{\iota \cup \gamma}  \arrow[swap]{d} & X \arrow{d}{f} \\
F_2 \times I  \arrow[swap]{r}{H \circ (k \times 1)} \arrow[dashed]{ur}{G} & Y \\
\end{tikzcd} \]

It can then be checked that the composite $F_1 \times I \to F_2 \times I \to X$, which defines $\gamma_1$, satisfies the required properties to show that $\gamma_1$ represents $\gamma$. $\hfill \square$

To figure out what this means for the actions of fundamental groups on homotopy fibres of an arbitrary composition of maps, we consider the following diagram:

\[ \begin{tikzcd}
X  \arrow{r}{f} & Y \arrow{r}{g} & Z \\
A \arrow{u} \arrow{d} \arrow{r}{i} & B \arrow{u} \arrow{d} \arrow{r}{j}  & C \arrow{u} \arrow{d} \\
A \times C^{I_+} \arrow{d} \arrow{r} & B \arrow{d} \arrow{r}  & C \arrow{d} \\
A \times B^{I_+} \times C^{I_+} \arrow[two heads]{r} & B \times C^{I_+} \arrow[two heads]{r} & C \\
\end{tikzcd} \]

where all vertical maps are weak equivalence, the map between the top two rows is a Reedy cofibrant approximation and the remaining maps between rows are the canonical ones. The bottom row gives a composite of fibrations which satisfies the conditions of the Lemma 2.2.3. Using the diagram, we can conclude that $\pi_1(Y)$ acts nilpotently on $H_n(Ff)$ iff  $\pi_1(B \times C^{I_+})$ acts nilpotently on $H_n(F_1)$. Similar conclusions hold for $Fg$ and $Fh$, where $h = gf$. We also have, with notation as in the previous lemma, that $\pi_1(F_2)$ acts trivially on $\pi_n(F_1)$ iff the induced map $Fh \to Fg$ induces a trivial action on its homotopy fibre, which is equivalent to $Ff$. 

\subsection{Properties of T-complete nilpotent groups}

Before moving on to the proof of the main theorem, we should say a few words about the category $\mathcal{B}_T$ of T-complete abelian groups. The definition we use states that an abelian, or nilpotent, group $G$ is T-complete if the T-completion map $G \to \mathbb{E}_T G$ is an isomorphism. Now, an abelian group $A$ is T-complete iff $\text{Hom}(\mathbb{Z}[T^{-1}],A) = \text{Ext}(\mathbb{Z}[T^{-1}], A) = 0$, see \cite[Prop. 10.1.18]{ref2},  and it follows from this that the kernel and cokernel of any homomorphism between T-complete abelian groups are T-complete. A nilpotent group is T-complete iff it is $\mathcal{B}_T$-nilpotent and so it follows that if $f: G \to H$ is a homomorphism between T-complete nilpotent groups, then the kernel and, if it exists, the cokernel of $f$ are T-complete, by \cite[Ch. III, Lemma 5.8]{ref9}. Any T-complete abelian group has the structure of a $\hat{\mathbb{Z}}_T$-module, and the universal property of completion implies that this $\hat{\mathbb{Z}}_T$-module structure is unique. Moreover, any homomorphism between T-complete abelian groups is automatically a $\hat{\mathbb{Z}}_T$-module homomorphism - in fact, it is a product of $\hat{\mathbb{Z}}_p$-module homomorphisms between the individual p-completions. Since each $\hat{\mathbb{Z}}_p$ is a PID, the kernel and cokernel of a homomorphism between $f \hat{\mathbb{Z}}_T$-modules are $f \hat{\mathbb{Z}}_T$-modules. As before, it follows that the kernel and cokernel, if it exists, of any homomorphism between $f \hat{\mathbb{Z}}_T$-nilpotent groups are $f \hat{\mathbb{Z}}_T$-nilpotent. Although we will not discuss it further in this paper, there is an evident alternative notion to finitely generated modules over $\hat{\mathbb{Z}}_T$, where we let an abelian group $A$ be in $\mathcal{F}_T$ iff $A$ is T-complete and, for every $p \in T$, $\hat{A}_p$ is an $f \hat{\mathbb{Z}}_p$-module. Again we have that the kernel and cokernel, if it exists, of a homomorphism between $\mathcal{F}_T$-nilpotent groups are $\mathcal{F}_T$-nilpotent. We record a few more such results in the following lemma:

\begin{Lemma}
	Let $G$ be a T-complete nilpotent group, and $H$ a T-complete subgroup. Then:
	
	i) there is a subnormal series $H = H_0 \leq H_1 \leq ... \leq H_k = G$, where each $H_i$ is T-complete, \\
	ii) if $G$ is f$\hat{\mathbb{Z}}_T$-nilpotent, then so is $H$, \\
	iii) if T is a finite set of primes and $G$ is $f \hat{\mathbb{Z}}_T$-nilpotent, then $G$ satisfies the ascending chain condition (ACC) for T-complete subgroups, \\
	iv) if $G$ is a T-torsion $f \hat{\mathbb{Z}}_T$-nilpotent group, then $G$ is finite.
	
\end{Lemma}

\textbf{Proof:} i) We will induct on the nilpotency class of $G$, noting that the result is trivial if $G$ is abelian. Let $e = G_0 \leq ... \leq G_q = G$ represent $G$ as a $\mathcal{B}_T$-nilpotent group. Let:

\[ \frac{H}{H \cap G_1} = K_0 \leq K_1 \leq ... \leq K_k = \frac{G}{G_1}
\]

be a subnormal series as is guaranteed to exist by the inductive hypothesis. Note, for example, that $H \cap G_1$ is T-complete since it is the kernel of $H \to \frac{G}{G_1}$. Let $\pi: G \to \frac{G}{G_1}$ denote the quotient, and define $H_{i+1} = \pi^{-1}(K_i)$. Then each $H_{i}$ is T-complete and $H_i$ is a normal subgroup of $H_{i+1}$. Moreover, $H$ is a normal subgroup of $H_1$, since $G_1$ is central in $G$. 

ii) If $\{G_i\}$ represents $G$ as an $f\hat{\mathbb{Z}}_T$-nilpotent group, then $\{H \cap G_i\}$ represents $H$ as an $f\hat{\mathbb{Z}}_T$-nilpotent group, since a T-complete submodule of an $f \hat{\mathbb{Z}}_T$-module is an $f \hat{\mathbb{Z}}_T$-module.

iii) When $G$ is abelian this follows from the fact that $\hat{\mathbb{Z}}_T$ is Noetherian. In general, this follows from a standard inductive argument on the nilpotency class of $G$.

iv) Since each $\hat{\mathbb{Z}}_p$ is $\{p\}$-local, we must have $\hat{G}_p = 1$ for all but finitely many primes. Therefore, we can reduce to the case where $T$ is a finite set of primes. When $G$ is abelian, the ACC implies that there is a product of primes in T, $r$, such that $rg = 0$ for all $g \in G$. Tensoring $\frac{\mathbb{Z}}{r\mathbb{Z}}$ with a suitable $\hat{\mathbb{Z}}_T$-free resolution of $G$, we conclude that $G$ is finite. The general case then follows by induction. $\hfill \square$

\newpage

\section{Completion preserves connected homotopy fibre squares}

\subsection{Proof of the main theorem}

The cocellular construction of the completion of a nilpotent space X can be modified in the following way. First, when replacing X by a Postnikov tower, we may assume the coattaching maps $X_i \to K(A,n)$ are cofibrations, by Lemma 2.2.2. Then to construct the completion we inductively use commutative squares of the form:

\[ \begin{tikzcd}
X_i \arrow[swap]{d}{} \arrow{r} &  \hat{X}_i \arrow{d}{} \\
K(A,n) \arrow{r} & \widehat{K(A,n)} \\
\end{tikzcd} \]

and the map $X_{i+1} \to \hat{X}_{i+1}$ is then defined as the canonical map between homotopy fibres.

\begin{Lemma} If $X$ is a connected nilpotent space and $f:X \to K(A,1)$ induces a surjection on $\pi_1$, then $A$ acts nilpotently on $H_{*}(Ff)$. \end{Lemma}

\textbf{Proof:} Since $f$ is surjective on $\pi_1$, there is a central series for $\pi_1(X)$ which realises $A$ as a quotient of $\pi_1(X)$ by a term of the central series. A simple diagram chase shows that we can take $f$ to be the projection onto a quotient tower, where $X$ is the limit of a Postnikov tower constructed using the chosen central series. We will show inductively that $A$ acts nilpotently on the homology of the fibres $F_i$ of the maps $d: X_i \to K(A,1)$. This can be shown by applying Lemma 2.2.3, and subsequent discussion, to the composites:

\[\begin{tikzcd}
X_{i+1} \arrow{d} \arrow{dr} \\
X_i \arrow{r} & K(A,1) \\
\end{tikzcd}  \]

Here, we assume that $A$ acts nilpotently on $H_{*}(F_i)$ and we know that $\pi_1(F_i)$ acts trivially on the homology of the fibre of $F_{i+1} \twoheadrightarrow F_i$, since this is a principal fibration and $\pi_1(F_{i+1}) \to \pi_1(F_i)$ is surjective. Therefore, an application of the Serre spectral sequence shows that $A$ acts nilpotently on $H_{*}(F_{i+1})$.  $\hfill \square$

\newpage

\begin{Proposition} Let $f:X \to A$ and $g: Y \to A$ be maps between connected nilpotent spaces such that $N(f,g)$ is connected. Let $\phi$ denote a completion as described above and let $\hat{\phi_T}$ denote an arbitrary completion. If we have a commutative diagram of the form:
	
	\[ \begin{tikzcd}
	X \arrow[swap]{d}{\hat{\phi}_T} \arrow{r}{f} & A \arrow{d}{\phi}  & Y \arrow[swap]{l}{g} \arrow{d}{\hat{\phi}_T}  \\
	\hat{X}_T \arrow[swap]{r}{\hat{f}} & \hat{A} & \hat{Y}_T \arrow{l}{\hat{g}}   \\
	\end{tikzcd} \] 
	
	then the induced map $N(f,g) \to N(\hat{f}, \hat{g})$ is a completion at T.
	
\end{Proposition}

\textbf{Proof:} We can assume that $A$ is the limit of a Postnikov tower with coattaching maps that are cofibrations as above. We define a filtration of $P := N(f,g)$ by pullbacks $P_i$:

\[ \begin{tikzcd}
P_i \arrow[swap]{d}{} \arrow{r} &  (A \times A) \times_{A_i \times A_i} A_i^{I_{+}}  \arrow{d}{} \\
X \times Y  \arrow{r} & A \times A \\
\end{tikzcd} \]

and let $Q_i$ denote the pullbacks in the corresponding filtration of $Q := N(\hat{f}, \hat{g})$. We have maps $P_i \to Q_i$ and we start by inductively showing that they are $\mathbb{F}_T$-equivalences. Since $P_0 = X \times Y$, the base case holds by assumption. Assume that $P_i \to Q_i$ is T-completion. We have a commutative square:

\[ \begin{tikzcd}
P_i \arrow[swap]{d}{} \arrow{r} &  Q_i \arrow{d}{} \\
\Omega K(B,m) \arrow{r} & \Omega \widehat{K(B,m)} \\
\end{tikzcd} \]

where we are using the induced tower structures on $A^{I_+}$ and $\hat{A}^{I_+}$ to construct the vertical maps. A diagram chase using our cocellular constructions of $A$ and $\hat{A}$ implies that the induced map between homotopy fibres is the map $P_{i+1} \to Q_{i+1}$. 

The Zeeman comparison theorem, or more specifically a refinement due to Hilton and Roitberg \cite{ref6}, can now be used either because $m > 2$ and so the base space is simply connected, or because $m = 2$ and  our previous work implies that $B$ and $\mathbb{E}_T B$ act nilpotently on the homology of the respective fibres. For the left hand fibration, this follows from Lemma 3.1.1 since $P$ connected implies that $P_i$ is connected and $\pi_1(P_i) \to B$ is surjective. 

For the right hand fibration, we are assuming that $P_i \to Q_i$ is T-completion. Therefore, $\pi_1(Q_i) \to \mathbb{E}_T B$ is surjective, since $\mathbb{E}_T$ is right exact. We can then apply Lemma 2.2.3 to the composite:

\[\begin{tikzcd}
Q_{i} \arrow[swap]{d}{u} \arrow{dr}{w} \\
\widehat{K(B,1)} \arrow[swap]{r}{v} & K(\mathbb{E}_T B,1) \\
\end{tikzcd}  \]

Lemma 3.1.1 tells us that $\mathbb{E}_T B$ acts nilpotently on the homology of $Fw$, and $\mathbb{E}_T B$ acts trivially on $Fv$, since $v$ is equivalent to a principal fibration. Therefore, an application of the Serre spectral sequence shows that $\mathbb{E}_T B$ acts nilpotently on the homology of $Q_{i+1}$. 

 Therefore, we can inductively show that $P \to Q$ is an $\mathbb{F}_T$-equivalence. We can also use this filtration of $Q$ to show that the homotopy groups of $Q$ are T-complete and so $Q$ is T-complete and $P \to Q$ is T-completion as desired. $\hfill \square$

We now consider functorial completions, as can be obtained via Bousfield localisation, for example. We have:

\begin{Theorem} Let $X \to \hat{X}_T$ denote a functorial completion. If a commutative square of connected nilpotent spaces:
	
	\[ \begin{tikzcd}
	P \arrow[swap]{d}{k} \arrow{r}{h} &  Y \arrow{d}{g} \\
	X \arrow[swap]{r}{f} & A\\
	\end{tikzcd} \]
	
	is a homotopy fibre square, then so is:
	
	\[ \begin{tikzcd}
	\hat{P}_T \arrow[swap]{d}{\hat{k}_T} \arrow{r}{\hat{h}_T} &  \hat{Y}_T \arrow{d}{\hat{g}_T} \\
	\hat{X}_T \arrow[swap]{r}{\hat{f}_T} & \hat{A}_T \\
	\end{tikzcd} \]	
\end{Theorem}

\textbf{Proof:} Using the functorial factorisation of any map as a cofibration followed by a weak equivalence, we may assume that each $X \to \hat{X}_T$ is a cofibration. We can also assume that each of the spaces $P,X,Y,A$ are $q$-cofibrant. We now put ourselves in the context of Proposition 3.1.2. Firstly, there is a weak equivalence $v: \hat{A}_T \to \hat{A}$ such that $v \hat{\phi}_T = \phi$ by the universal property of completion and the assumption that $\hat{\phi}_T$ is a cofibration. We now let $\hat{f} = v \hat{f}_T, \hat{g} = v \hat{g}_T$ and these choices make the diagram in the statement of Proposition 3.1.2 commute. Therefore, $N(f,g) \to N(\hat{f},\hat{g})$ is T-completion. It follows easily that $N(f,g) \to N(\hat{f}_T, \hat{g}_T)$ is also T-completion. It now follows that $\hat{P}_T \to N(\hat{f}_T, \hat{g}_T)$ is an $\mathbb{F}_T$-equivalence and, therefore, a weak equivalence as desired. $\hfill \square$

\newpage

\subsection{Fracture Theorem as a Consequence}

In order to deduce the fracture theorem as a corollary of our results, we first need to show that the homotopy pullback in question is connected. This is the content of the following lemma:

\begin{Lemma}
	If $G$ is a $T$-local nilpotent group, then the function $\varphi: \mathbb{E}_T G \times G_0 \to (\mathbb{E}_T G)_0$ is surjective.
\end{Lemma}

\textbf{Proof:} We first assume that $G$ is abelian. Let $J$ denote the image of $\varphi$, and $P$ denote the kernel of $\varphi$, which is just the evident pullback. We will first show that $J$ is rational. Since $J$ is a subgroup of $(\mathbb{E}_T G)_0$, we have $\mathbb{H}_T J = 0$. Therefore, we have a short exact sequence:

\[ 0 \to \mathbb{E}_T P  \to \mathbb{E}_T G \to \mathbb{E}_T J \to 0
\]

The universal property of the pullback implies that the first map splits and, in particular, is surjective. It follows that $J$ is T-local and $\mathbb{E}_T J = \mathbb{H}_T J = 0$. Therefore, $J$ is rational. Now  $\varphi_0$ is surjective and factors through $J$, so $J = (\mathbb{E}_T G)_0$ and $\varphi$ is surjective.

The result for general nilpotent groups $G$ can now be proven inductively on the nilpotency class of $G$, using \cite[Lemma 7.6.1]{ref2}. A key point is that the image of $\mathbb{E}_T Z(G)$ in $\mathbb{E}_T G$ is a central subgroup - this can be seen from the Postnikov tower construction of completion applied to the upper central series of $G$. $\hfill \square$

Finally, we give the proof of the fracture theorem that we have been building toward:

\begin{Theorem}
	Let $X$ be a $T$-local connected nilpotent space. Then any commutative square:
	
\[ \begin{tikzcd}
X \arrow{r}{\hat{\phi}} \arrow[swap]{d}{\psi} & \hat{X}_T \arrow{d}{\phi} \\
X_0 \arrow[swap]{r}{\hat{\phi}_0} & (\hat{X}_T)_0 \\
\end{tikzcd}
\]

with $\hat{\phi}$ a T-completion and $\psi,\phi$ rationalisations, is a homotopy fibre square.

\end{Theorem}

\textbf{Proof:} We can assume that $\phi$ is a fibration. Then, by Lemma 3.2.1, the pullback $P$ is connected and we have a comparison map $f: X \to P$. By Theorem 3.1.3, applying functorial rationalisation and completion shows that $f_0$ and $\hat{f}_T$ are weak equivalences, respectively. It follows that $\tilde{H}_{*}(Cf; \mathbb{Q}) = \tilde{H}_{*} (Cf; \bigoplus_{p \in T} \mathbb{F}_p) = 0$. Now $\tilde{H}_{*}(Cf)$ is T-local, since $X$ and $P$ are, and it follows that $\tilde{H}_{*}(Cf)$ is also local away from T and has trivial rationalisation. Therefore, $\tilde{H}_{*}(Cf) = 0$ and $f$ is a homology isomorphism between connected nilpotent spaces, so must be a weak equivalence. $\hfill \square$

\printbibliography

\end{document}